\documentclass[12pt]{article}
\usepackage{amsfonts}
\usepackage{amsmath,amsthm}
\usepackage[all,frame]{xy}
\usepackage{enumerate}
\usepackage{multirow}
\usepackage{geometry}
\numberwithin{equation}{section}

\newtheorem{thm}{Theorem}[section]

\newtheorem{proposition}[thm]{Proposition}

\theoremstyle{definition}

\newenvironment{prf}[1][Proof]{\begin{trivlist}
\item[\hskip \labelsep {\bfseries #1}]}{\end{trivlist}}
\newenvironment{remark}[1][Remark]{\begin{trivlist}
\item[\hskip \labelsep {\bfseries #1}]}{\end{trivlist}}

\newcommand{\qd}{\ensuremath{\partial}}
\newcommand{\zh}{\Theta_{V}}
\newcommand{\wtil}{\ensuremath{\tilde{\omega}}}
\newcommand{\der}[1]{\ensuremath{\partial_{#1}}}
\newcommand{\ub}[2]{\ensuremath{_{#1}^{(#2)}}}
\newcommand{\mc}[1]{\ensuremath{\mathcal{#1}}}

\newcommand{\vv}{\mathbf{1}}
\newcommand{\Tr}{\mathrm{Tr}\, }
\newcommand{\perm}{\mathrm{perm}}
\newcommand{\pperm}{\mathrm{pperm}}
\newcommand{\dom}{\mathrm{dom}\,}

\begin{document}

\title{Virasoro Correlation Functions for Vertex Operator Algebras}
\author{
Donny Hurley\thanks{
Supported by the Science Foundation Ireland Research Frontiers Programme} 
\ and Michael P. Tuite \\
School of Mathematics, Statistics and Applied Mathematics \\
National University of Ireland Galway, \\
University Road,
Galway, Ireland.}
\maketitle

\begin{abstract}
We consider all genus zero and genus one correlation functions for the  Virasoro vacuum descendants of a vertex operator algebra. 
These are described in terms of explicit generating functions
that can be combinatorially expressed in terms of graph theory related to derangements in the genus zero case 
and to partial permutations in the genus one case.
\end{abstract}

\newpage

\section{Introduction}
\label{Intro}
A Vertex Operator Algebra (VOA) (e.g. \cite{B}, \cite{FLM}, \cite{K})  is an algebraic structure closely related to Conformal Field Theory (CFT) in physics e.g. \cite{BPZ}, \cite{DMS}. 
An essential ingredient of a VOA or CFT is the existence of a Virasoro vector whose vertex operator modes generate a Virasoro subalgebra of central charge $C$.  
The purpose of this paper is to describe all genus zero and genus one correlation functions for Virasoro descendants of the vacuum vector in terms of explicit generating functions.
Examples of such Virasoro correlation functions have long been known in the VOA and CFT literature.
For instance, they occur in describing the Kac determinant, CFT Ward identities and VOA modular differential equations e.g. \cite{DMS}, \cite{Z}, \cite{MT2}.
Our aim here is to provide a complete description of these Virasoro correlation functions at genus zero and one.   

We begin in Section~\ref{genus0sect} with a brief review of the theory of Vertex Operator Algebras (VOA). 
We describe a set of symmetric rational generating functions for all genus zero correlation functions for Virasoro vacuum descendants.
These satisfy a recursion relation due to Zhu \cite{Z} which in Theorem~\ref{g0mainthm} imply that
each generating function can be expressed as the sum of certain rational weights of appropriate graphs related to derangements in combinatorics. 
We also show that the generating function can be combinatorially expressed in terms of a $\beta$-extension of a permanent \cite{FZ}.  

In Section~\ref{genus1chap} we consider genus one correlation functions for Virasoro vacuum descendants.
We begin with a review of elliptic  and modular functions and prove in Theorem~\ref{g1breakformula} that the Weierstrass 
function satisfies an important non-linear differential equation.
(As a by-product of this analysis we also obtain a closed formula for 
the modular derivative of any Eisenstein series).
We then describe a set of symmetric elliptic generating functions for all such genus one Virasoro correlation functions which satisfy a genus one Zhu recursion relation \cite{Z}.
The resulting expression for the generating function is not manifestly symmetric. However, by use of the non-linear differential equation of 
Theorem~\ref{g1breakformula} we show in 
Theorem~\ref{g1_maintheorem} how to express each generating function as a sum of weights of appropriate graphs related to partial permutations.
We also show that this generating function can be expressed in terms of an $\alpha,\beta$-extension of a partial permanent 
defined in Section~\ref{sec:pperm}, a natural generalization to partial permutations of the notion of the 
$\beta$-extension of a permanent defined for permutations. 
The Section concludes with an alternative graphical approach based on permutations.
Finally, the Appendix describes the computation of the exponential generating function for enumerating the various graphs employed in the earlier sections.

\medskip

\section{Genus Zero Virasoro Correlation Functions}
\label{genus0sect}
\subsection{ Vertex operator algebras}
\label{VOA}
We review some basic properties of Vertex Operator Algebras  e.g. \cite{FLM}, \cite{FHL}, \cite{K}, \cite{LL}, \cite{MN}, \cite{MT2}.
A Vertex Operator Algebra (VOA) is a quadruple $(V, Y, \vv, \omega)$
consisting of a $\mathbb{Z}$-graded complex vector space $V =\bigoplus_{n\in\mathbb{Z}}V_n$ with $\dim V_{n}<\infty$, a linear map $Y \rightarrow (\mathrm{End}V)[[z, z^{-1}]$, for formal parameter $z$ and a pair of states: the vacuum $\mathbf{1} \in V_0$ and a conformal vector $\omega \in V_2$. For each $v \in V$ we have a vertex operator
\begin{equation*}
Y(v,z) = \sum_{n\in\mathbb{Z}}v_n z^{-n-1},
\end{equation*}
which satisfies the following axioms: 
\begin{description}
\item[Lower Truncation.]
For each pair $u,v\in V$ then $v_{n}u=0$ for $n\gg 0$.
\item[Locality.] 
$(z_1-z_2)^N[Y(u,z_1),Y(v,z_2)] = 0$ for $N\gg0$.
\item[Creativity.] 
$Y(v,z)\vv = v + O(z)$.
\item[Virasoro Structure.] 
There exists a conformal vector $\omega\in V_2$ with
\begin{displaymath}
Y(\omega,z) = \sum_{n\in\mathbb{Z}}L_n z^{-n-2},
\end{displaymath}
where the modes $L_{n}=\omega_{n+1}$ satisfy a Virasoro Algebra of central charge $C$ with bracket relation
\begin{align*}
[L_m,L_n] =(m - n)L_{m+n} + 
\frac{(m^3 - m)}{12}C\delta_{m,-n}.
\end{align*} 
Furthermore, the $\mathbb{Z}$-grading on $V$ is provided by $L_{0}$ with $L_{0}v=nv$ for all $v\in V_{n}$.
\item[Translation.] 
$[L_{-1}, Y(v, z)] = \partial_z Y(v, z)$.
\end{description} 
\medskip

\noindent These axioms imply $Y(\vv,z)=\mathrm{Id}_{V}$ and associativity:
\begin{eqnarray}
(z_{0}+z_{2})^M Y(u,z_{0}+z_{2})Y(v,z_{2})w &=& (z_{0}+z_{2})^M Y(Y(u,z_{0})v,z_{2})w,
\label{Assoc}
\end{eqnarray}
for $u,v,w\in V$ and $M \gg 0$  \cite{FHL}. In this paper we consider VOAs of CFT-type for which $V_{0}=\mathbb{C}\vv$ (and consequently, $V_{n}=0$ for $n<0$ \cite{MT2}). 

\subsection{Virasoro correlation functions}\label{sec:g0vircorrfn}
Let  $v_{1},\ldots,v_{n}\in V$ and
define the genus zero $n$-point correlation function for formal
parameters $z_{1},\ldots,z_{n}$ by
\begin{equation}
Z_V^{(0)}\left(v_{1},z_1;\ldots;v_{n},z_n\right)=\langle 
\vv^\prime, Y(v_{1},z_{1})\ldots Y(v_{n},z_{n})\vv
\rangle,
\label{eq:npt-Fun}
\end{equation}
where $\vv^\prime$ denotes the dual of the vacuum vector.
Lower truncation and locality imply that the $n$-point function is a rational function in $z_{1},\ldots,z_{n}$ (expanded in the formal domain
$|z_{1}|>\ldots>|z_{n}|$). Thus 
it may be considered as a rational function on the complex Riemann sphere for $z_{i}\in \hat{\mathbb{C}}$ \cite{FHL}, \cite{MT2}.     
In particular, we consider the genus zero Virasoro $m$-point function 
\begin{align}\label{g1_gndefn}
G_m^{(0)}(z_1,\ldots,z_m)=
Z_V^{(0)}\left(\omega,z_1;,\ldots;\omega,z_m\right),
\end{align}
with $v_{1}=\ldots=v_{m}=\omega$, the Virasoro vector. We then find
\begin{proposition}
\label{prop_Ggen0}
$G_m^{(0)}(z_1,\ldots,z_m)$ is a symmetric rational function in $z_1,\ldots,z_m$. The set of such functions generate all genus zero $n$-point correlation functions for Virasoro vacuum descendants.
\end{proposition}
\begin{remark} 
We refer to $G_m^{(0)}(z_1,\ldots,z_m)$ as the order $m$ genus zero Virasoro generating function.
\end{remark}
\begin{prf}
Locality implies that $G_m^{(0)}(z_1,\ldots,z_m)$ is a 
symmetric rational function in $z_1,\ldots,z_m$.  
We illustrate the generating property by considering the 2-point function for two Virasoro vacuum descendants $v_1=L_{-l_1}\ldots L_{-l_r}\vv$ and $v_2=L_{-k_1}\ldots L_{-k_s}\vv$. 
Consider 
\begin{equation}\label{g0_heq1}
\langle
\vv^\prime,
Y(Y(\omega,w_1)\ldots Y(\omega,w_r)\vv,w)Y(Y(\omega,z_1)\ldots Y(\omega,z_s)\vv,z)\vv\rangle.
\end{equation}
The two point function for $v_1,v_2$ is the coefficient of $\prod_{i=1}^{r}w_{i}^{l_{i}-1}\prod_{j=1}^{s}z_{j}^{k_{j}-1}$ in \eqref{g0_heq1}. By repeatedly using associativity \eqref{Assoc} we find
\begin{eqnarray*}
&&
\prod_{i=1}^{r}(w_{i}+w)^{M}\prod_{j=1}^{sn}(z_{j}+z)^{N}
\langle\vv^\prime,
Y(Y(\omega,w_{1})\ldots Y(\omega,w_r)\vv,w)Y(Y(\omega,z_1)\ldots Y(\omega,z_s)\vv,z)\vv\rangle
\\
&&
=\prod_{i=1}^{r}(w_{i}+w)^{M}\prod_{j=1}^{s}(z_{j}+z)^{N}
\langle\vv^\prime,
Y(\omega,w_1+w)Y(Y(\omega,w_2)\ldots Y(\omega,w_m)\vv,w)
\\
&&
\quad Y(\omega,z_1+z)Y(Y(\omega,z_2)\ldots Y(\omega,z_n)\vv,z)\vv\rangle
\\
&&
=
\langle\vv^\prime,
\prod_{i=1}^{r}(w_{i}+w)^{M}\prod_{j=1}^{s}(z_{j}+z)^{N}
Y(\omega,w_{i}+w)Y(\omega,z_{j}+z)\vv\rangle
\\
&&
=\prod_{i=1}^{r}(w_{i}+w)^{M}\prod_{j=1}^{s}(z_{j}+z)^{N}
G_{r+s}^{(0)}(w_1+w,\ldots,w_r+w,z_1+z,\ldots,z_s+z),
\end{eqnarray*}
for $M,N\gg 0$. 
Thus $\eqref{g0_heq1}$ and $G_{r+s}^{(0)}(w_1+w,\ldots,w_r+w,z_1+z,\ldots,z_s+z)$ are determined by the same rational function whose coefficient of $\prod_{i=1}^{r}w_{i}^{l_{i}-1}\prod_{j=1}^{s}z_{j}^{k_{j}-1}$  is  the 2-point function for $v_1,v_2$. 
The general result follows from a similar argument. \qed
\end{prf}

By Lemma~2.2.1 of \cite{Z} we obtain a recursion formula for the order $n$ Virasoro generating function: 
\begin{proposition}
\label{zhuredform}
\begin{align}\notag
G_n^{(0)}(z_1,\ldots,z_n) = 
\sum_{k=2}^n &\left(\frac{1}{z_{1k}}\der{z_k}
+\frac{2}{z_{1k}^2}\right)G_{n-1}^{(0)}(z_2,\ldots,z_n)\\ &
+ \frac{C}{2}\frac{1}{z_{1k}^4}G_{n-2}^{(0)}(z_2,\ldots,\hat{z_k},\ldots,z_n),
\label{g0reductionformula} 
\end{align}
where $z_{ij}=z_i - z_j$ and $\hat{z_k}$ denotes that $z_k$ is omitted. \qed
\end{proposition}
Thus we find using $G_0^{(0)}=1$ and $G_1^{(0)}(z_1)=0$ that
\begin{equation*}
G_2^{(0)}(z_1,z_2) = \frac{C}{2}\frac{1}{z_{12}^4},\quad G_3^{(0)}(z_1,z_2,z_3) = \frac{C}{z_{12}^2z_{13}^2z_{23}^2}.
\end{equation*}

\subsection{Genus zero Virasoro graphs}
\label{ssect:genus0}
We next develop a graphical/combinatorial approach for computing the order $n$ Virasoro generating function. 
To this end we define an \textit{Order $n$ Genus Zero Virasoro Graph}
to be a directed graph $\mathcal{G}$ with $n$ vertices labelled by $z_1,\ldots,z_n$ with the following properties.  
Each vertex has degree 2 with unit indegree and outdegree and
no vertex is joined to itself.
Thus the connected subgraphs of $\mathcal{G}$ are $r$-cycles with $r\ge 2$ vertices.
We define a weight $W(\mathcal{G})$ on $\mathcal{G}$ as follows. For each directed edge 
$\mathcal{E}_{ij}$, with $i \neq j$, 
\begin{center}$\xy(0,0)*{\cir<3pt>{}}="a"*+!R{z_i\,}; (15,0)*{\cir<3pt>{}}="b"*+!L{\,z_j}; \ar "b";"a";\endxy$
\end{center} 
define an edge weight
\begin{equation*}
W(\mathcal{E}_{ij})=\frac{1}{z_{ij}^2}.
\end{equation*}
Let $K$ be the number of cycles in \mc{G} and define 
\begin{equation*}
W(\mathcal{G})=\left(\frac{C}{2}\right)^K\prod_{\mathcal{E}}W(\mathcal{E}),
\end{equation*}
where the product ranges over all the edges of \mc{G}.

In the Appendix we describe a generating function for counting the number
$d_{n,K}$ of inequivalent order $n$ graphs with $K$ cycles. 
This is determined by the polynomial
\begin{equation}
d_{n}(\beta)=\sum_{K\ge 1}d_{n,K}\beta^K=(-1)^nn! \sum_{i=0}^{n} \frac{\beta^i}{i!}\binom{-\beta}{n-i},
\label{dnbeta}
\end{equation}
for cycle counting parameter $\beta$. 
Thus for $n=2$ we find $d_{2}(\beta)=\beta$ so that there is one Virasoro graph \mc{G} with weight 
\begin{equation*}
W(\mathcal{G})=\frac{C}{2}.\frac{1}{z_{12}^2}.\frac{1}{z_{12}^2}= G_2^{(0)}(z_1,z_2),
\end{equation*}
the order 2 Virasoro generating function. For
$n=3$ we have $d_3(\beta) = 2\beta$ so that there are two order $3$ Virasoro graphs:  
\begin{equation*}
\mathcal{G}_1=\xy (0,-2)*[o]=<0.5pt>+{\cir<3pt>{}}="a"*+!R{\,z_3}; (10,-2)*[o]=<0.5pt>+{\cir<3pt>{}}="b"*+!L{z_2\,}; (5,5)*[o]=<0.5pt>+{\cir<3pt>{}}="c"*+!D{z_1};
\ar "c";"a";\ar "b";"c";\ar "a";"b";
\endxy,
\qquad 
\mathcal{G}_2=\xy (0,-2)*[o]=<0.5pt>+{\cir<3pt>{}}="a"*+!R{\,z_2}; (10,-2)*[o]=<0.5pt>+{\cir<3pt>{}}="b"*+!L{z_3\,}; (5,5)*[o]=<0.5pt>+{\cir<3pt>{}}="c"*+!D{z_1};
\ar "c";"a";\ar "b";"c";\ar "a";"b";
\endxy.
\end{equation*}
The sum of the weights of these two graphs is
\begin{align*}
W(\mathcal{G}_1)+W(\mathcal{G}_2)&=\frac{C}{2}\frac{1}{z_{12}^2z_{23}^2z_{31}^2} + \frac{C}{2}\frac{1}{z_{13}^2z_{32}^2z_{21}^2}
= G_3^{(0)}(z_1,z_2,z_3),
\end{align*}
the order 3 Virasoro generating function. For $n=4$ one finds
$d_4(\beta) = 3\beta^2+6\beta$ so that there are 9 independent graph configurations whose weights sum to the order 4 Virasoro generating function.
These examples illustrate the following general result:
\begin{thm}\label{g0mainthm}
The  order $n$ genus zero Virasoro generating function $G_n^{(0)}$ is the sum of the weights of all inequivalent order $n$ genus zero Virasoro graphs. 
\end{thm}

\begin{prf}
We prove the result by induction in $n$. The result holds for $n=2$ and $3$ as described above.
For $n\ge 4$ we reinterpret the recursion formula equation~\eqref{g0reductionformula} in terms of a construction
of order $n$ Virasoro graphs from order $n-1$ and $n-2$ Virasoro graphs.

Divide the set of inequivalent order $n$ Virasoro graphs into the following two types:
\begin{enumerate}[Type I.]
\item Graphs for which the $z_1$ vertex belongs to a 2-cycle.
\item Graphs for which the $z_1$ vertex belongs to an $r$-cycle for $r\ge 3$.
\end{enumerate}
Consider the following term arising in the RHS of  \eqref{g0reductionformula} for $k=2,\ldots,n$
\begin{equation}
\label{pfg0eq2}
\frac{C}{2}\frac{1}{z_{1k}^4}
G_{n-2}^{(0)}(z_2,\ldots,\hat{z_k},\ldots,z_n).
\end{equation}
By induction, $G_{n-2}^{(0)}$ is the sum of weights of inequivalent Virasoro order $n-2$ graphs 
$\{\mathcal{G}^{n-2}\}$ with vertices labelled by $z_2,\ldots,\hat{z_k},\ldots,z_n$.
For each such $\mathcal{G}^{n-2}$, let $\mathcal{G}_{\mathrm{I}}^{n}$ denote the order $n$ Virasoro graph formed by the product of $\mathcal{G}^{n-2}$ and the 2-cycle with $z_1$ and $z_k$--vertices. 
Clearly, $\mathcal{G}_{\mathrm{I}}^{n}$ is an order $n$ graph of Type~I and all 
inequivalent such graphs arise in this way.   
Since
\begin{equation*}
W(\mathcal{G}_{\mathrm{I}}^{n})=\frac{C}{2}\frac{1}{z_{1k}^4}W(\mathcal{G}^{n-2}),
\end{equation*}
then summing \eqref{pfg0eq2} over all $k$ results in 
\begin{align}
\sum_{k\ge 2}\frac{C}{2}\frac{1}{z_{1k}^4}G_{n-2}^{(0)}(z_2,\ldots,\hat{z_k},\ldots,z_n) = \sum_{\mathcal{G}_{\mathrm{I}}^{n}} W(\mathcal{G}_{\mathrm{I}}^{n}),\label{g0_pftype1gph}
\end{align}
where $\mathcal{G}_{\mathrm{I}}^{n}$ ranges over the inequivalent Virasoro order $n$ graphs of Type I. 

\bigskip
We next consider the remaining terms arising on the RHS of  \eqref{g0reductionformula}
for $k=2,\ldots,n$:
\begin{equation*}
\mathcal{D}G_{n-1}^{(0)}(z_2,\ldots,z_n),
\end{equation*}
where $
\mathcal{D}=\sum_{k=2}^n(\frac{1}{z_{1k}} \der{ z_k}
+ \frac{2}{z_{1k}^2})\label{pfg0eq1}$.
By induction, $G_{n-1}^{(0)}(z_2,\ldots,z_n)$ is the sum of weights of all inequivalent order $n-1$ Virasoro graphs $\mathcal{G}^{n-1}$ labelled by $z_2,\ldots,z_n$. 
Let $\{\mathcal{E}_{ij}\}$ denote the edges for such an order $n-1$ graph. 
Recalling that $W(\mathcal{E}_{ij})=\frac{1}{z_{ij}^2}$ 
we find that
\begin{equation*}
\mathcal{D}W(\mathcal{G}^{n-1}) = 
\sum_{\mathcal{E}_{ij}}W(\mathcal{G}^{n-1})z_{ij}^2
\left(\frac{1}{z_{1i}}\der{ z_i}+ \frac{1}{z_{1j}}\der{ z_j} + \frac{1}{z_{1i}^2} + \frac{1}{z_{1j}^2}\right)\frac{1}{z_{ij}^2}.
\end{equation*}
Using the differential equation 
\begin{equation}\label{g0breakformula}
\left(-\frac{1}{x}\der{x} - \frac{1}{y}\der{y} + \frac{1}{x^2} + \frac{1}{y^2}\right)\frac{1}{(x-y)^2} = \frac{1}{x^2y^2},
\end{equation}
for $x=z_{1i}$ and $y=z_{1j}$ we obtain
\begin{align*}
\mathcal{D}W(\mathcal{G}^{n-1}) &= \sum_{\mathcal{E}_{ij}}W(\mathcal{G}^{n-1})\frac{z_{ij}^2}{z_{1i}^2 z_{1j}^2}
= \sum_{\mathcal{G}_{ij}^{n}}W(\mathcal{G}_{ij}^{n}),
\end{align*}
where $\mathcal{G}_{ij}^{n}$ denotes the order $n$ Virasoro  graph obtained by ``inserting'' a vertex labelled $z_1$ between each connected pair $z_i$ and $z_j$ in $\mathcal{G}^{n-1}$ i.e. schematically
\begin{equation*}
\cdots
\xy  
(0,0)*[o]=<0.5pt>+{\cir<3pt>{}}="a"*+!R{z_i\,}; (15,0)*[o]=<0.5pt>+{\cir<3pt>{}}="b"*+!L{\,z_j};
\ar "b";"a";
\endxy
\cdots
\quad
 \longrightarrow 
\quad\cdots
 \xy 
 (0,0)*[o]=<0.5pt>+{\cir<3pt>{}}="a"*+!R{z_i\,}; (15,0)*[o]=<0.5pt>+{\cir<3pt>{}}="e"*+!D{z_1}; (30,0)*[o]=<0.5pt>+{\cir<3pt>{}}="b"*+!L{\,z_j};
\ar "e";"a";
\ar "b";"e";
\endxy
\cdots
\end{equation*} 
Clearly, each order $n$ Virasoro graph $\mathcal{G}_{ij}^{n}$ is of Type~II.
Summing over all inequivalent order $n-1$ graphs  $\{\mathcal{G}^{n-1}\}$ results in a sum over all inequivalent graphs  $\lbrace\mathcal{G}_{\mathrm{II}}^{n}\rbrace$
of Type~II i.e.
\begin{align*}
\mathcal{D}G_{n-1}^{(0)}(z_2,\ldots,z_n) = \sum_{\mathcal{G}_{\mathrm{II}}^{n}}W(\mathcal{G}_{\mathrm{II}}^{n}).
\end{align*}
Combined with the Type~I graphs from \eqref{g0_pftype1gph} we therefore find
\begin{align*}
G_n^{(0)}(z_1,\ldots,z_n) = \sum_{\mathcal{G}_{\mathrm{I}}^{n}}W(\mathcal{G}_{\mathrm{I}}^{n}) + \sum_{\mathcal{G}_{\mathrm{II}}^{n}}W(\mathcal{G}_{\mathrm{II}}^{n}),
\end{align*}
as required.\qed
\end{prf}

\subsection{Derangements and the $\beta$-extension of a permanent}\label{sec:g0combin}
The set of  order $n$ Virasoro graphs $\lbrace\mathcal{G}\rbrace$ is in 1-1 correspondence with the set of derangements (i.e. fixed point free permutations) $F(\Phi)$ of the label set $\Phi=\lbrace 1,\ldots,n\rbrace$. 
Thus the directed edges $\{\mathcal{E}_{ij}\}$ of $\mathcal{G}$ correspond to a derangement permutation $\phi\in F(\Phi)$ with $\phi(i)=j\neq i$. In particular, connected cycle subgraphs correspond to cycles of the corresponding permutation.
Hence we may write the result of Theorem~\ref{g0mainthm} as 
\begin{equation}
G_n^{(0)}(z_1,\ldots,z_n) =\sum_{\phi\in F(\Phi)}\left(\frac{C}{2}\right)^{K(\phi)}\prod_{i=1}^n\frac{1}{z_{i\phi(k_i)}^2},
\label{g0_thmderange}
\end{equation}
where $K(\phi)$ is the number of cycles in the derangement permutation $\phi$.
The expression on the RHS of \eqref{g0_thmderange} can be expressed in terms the so-called \emph{$\beta$-extension} of a permanent of an $n\times n$ matrix $A = (a_{ij})$ defined by \cite{FZ}
\begin{equation}
\perm_\beta\,A = \sum_{\pi\in \Sigma_{n}}\beta^{K(\pi)}\prod_{i=1}^n a_{i\pi(i)},\label{g0_betaperm}
\end{equation}
where $\pi$ ranges over all permutations $\Sigma_{n}$ of the set $\lbrace1,\ldots,n\rbrace$ and $K(\pi)$ is the number of cycles in $\pi$. 
(For $\beta=1$ then $\perm_{1}\,A=\perm\,A$, the standard permanent, 
whereas for $\beta=-1$ we obtain $\perm_{-1}\,A=(-1)^n\det A$).
 In the present context, we define a symmetric matrix $A$ with components
\begin{align*}
a_{ij}= \left\{
  \begin{array}{ll}
    0 & i=j,\\
   \frac{1}{z_{ij}^2} & i\neq j.
  \end{array}\right.
\end{align*}
Since the diagonal elements of $A$ are $0$, only fixed point free permutations contribute to \eqref{g0_betaperm} so that Theorem~\ref{g0mainthm} can be restated as
\begin{thm}\label{g0mainthmperm}
The order $n$ genus zero Virasoro generating function is given by
\begin{displaymath}
G_n^{(0)}(z_1,\ldots,z_n) = \perm_{\frac{C}{2}}\,A.\qed
\end{displaymath}
\end{thm}
\bigskip

\section{Genus One Virasoro Correlation Functions}
\label{genus1chap}
\subsection{Elliptic functions and modular forms}\label{Genus1Elliptic}
We first consider aspects of the theory of elliptic and modular functions relevant to our discussion. We note the following conventions: $q_z \equiv e^{z}$, $q \equiv q_{2\pi i\tau} =e^{2\pi i \tau}$ and $\qd \equiv q\der{q} = \frac{1}{2\pi i}\der{\tau}$ for $z\in \mathbb{C}$ and $\tau\in\mathbb{H}$, the complex upper half plane.
Define 
\begin{eqnarray}
P_{2}(z,q) &=&\wp ( z,\tau)+E_{2}(q )  \notag \\
&=&\frac{1}{z^{2}}+\sum_{k=2}^{\infty }(k-1)E_{k}(q )z^{k-2},  \label{P2}
\end{eqnarray}%
where $\wp (z,\tau)$ is the Weierstrass function periodic in $2\pi i$ and $2\pi i\tau$ e.g. \cite{WW}, \cite{C}. $E_{k}(q)$ is zero for $k$ odd, and for $k$ even is the Eisenstein series 
\begin{equation}
E_{k}(q)=-\frac{B_{k}}{k!}+\frac{2}{(k-1)!}\sum\limits_{n\geq 1}\frac{n^{k-1}q^{n}}{1-q^{n}},  \label{intro_Ekdefn}
\end{equation}%
where $B_k$ is the Bernoulli number determined by
\begin{equation*}
\frac{z}{q_{z}-1}= \sum\limits_{k\geq 0} \frac{B_{k}}{k!}z^k = 1-\frac{1}{2}z+\frac{1}{12}z^2+ \hdots
\end{equation*}
$E_k(q)$ is a holomorphic modular form of weight $k$ on $SL(2, \mathbb{Z})$ for $k\ge 4$ \cite{S} whereas $E_2(q)$ is a quasi-modular form \cite{KZ}.
We also define $P_k(z,q)$ for $k\ge 0$ by
\begin{align}
P_0(z,q) &= -\log(z) + \sum_{k\ge 2} \frac{1}{k}E_k(q)z^k,\\
P_k(z,q) &= -\frac{1}{k-1}\der{z}P_{k-1}(z,q)
= \frac{1}{z^k} + (-1)^k\sum_{n\ge k}^{\infty}E_n(q)\binom{n-1}{k-1}z^{n-k}.\label{intro_Pkdefn}
\end{align}
$P_1$ is periodic in $2\pi i$ but is not periodic in $2\pi i\tau$ but $P_k$ is periodic in both $2\pi i$ and $2\pi i \tau$ for all $k\ge 2$.  
$P_0(z)$ is related to the elliptic prime form $K(z,q)$ \cite{Mu}
\begin{align}
K(z,q) = \exp(-P_0(z,q)) = -\frac{i\theta_1(z,q)}{\eta^3(q)},
\label{g1_Kthetaeqn}
\end{align}
where 
\begin{eqnarray}
\theta_1(z,q) &=& i\sum_{n\in\mathbb{Z}}q^{\frac{1}{2}(n+\frac{1}{2})^2}q_z^{n+\frac{1}{2}}(-1)^n,
\label{intro_thetadef}\\
\eta(q) &=& q^{1/24}\prod_{n \geq 1}(1-q^n),
\label{intro_dedekindef}
\end{eqnarray}
are the Jacobi theta series and the Dedekind $\eta$-function respectively.

We find the following:
\begin{proposition}\label{genus1ellthm}
\begin{equation}
2\qd P_0(z,q)=P_2(z,q)-P^2_1(z,q)-3E_2(q).
\label{genus1ellthmeqn}
\end{equation}
\end{proposition}
\begin{prf}
We first note from \eqref{intro_Ekdefn} and \eqref{intro_dedekindef} that 
\begin{equation}
\qd\log(\eta(q)) = -\frac{1}{2} E_2(q).
\label{deleta}
\end{equation}  
Furthermore, $\theta_1(z,q)$ obeys the heat equation 
\begin{equation}
2\qd\theta_1(z,q) = \der{z}^2\theta_1(z,q).
\label{g1_theta1_3}
\end{equation} 
It thus follows from \eqref{g1_Kthetaeqn} that
\begin{eqnarray*}
2\qd\theta_1(z,q) &=& \left(-2\qd P_0(z,q) - 3E_2(q)\right)\theta_1(z,q),\\
\der{z}^2\theta_1(z,q) &=& \left(P_1^2(z,q)-P_2(z,q)\right)\theta_1(z,q).
\end{eqnarray*}
Hence the heat equation \eqref{g1_theta1_3} implies the result. \qed 
\end{prf}
Comparing the Laurent expansions in $z$ of both sides of \eqref{genus1ellthmeqn} immediately implies
\begin{proposition}\label{g1_E2cor} 
The Eisenstein series obeys
\begin{align}
\qd E_k(q) = \frac{k(k+3)}{2}E_{k+2}(q) - \frac{k}{2}\sum_{r\ge 2}^kE_r(q)E_{k+2-r}(q)
\label{g1_E2pfeqn}.\qed
\end{align}
\end{proposition}
\begin{remark}
Let $D_{k}=\qd  + kE_{2}(q)$ denote the modular derivative 
which when applied to a modular form $f(\tau)$ of weight $k$ results
in a modular form $D_kf(\tau)$ of weight $k+2$ \cite{Z}, \cite{Ma}. 
For $k\ge 4$, Proposition \ref{g1_E2cor} implies  
\begin{align*}
D_{k}E_k(q) = \frac{k(k+3)}{2}E_{k+2}(q) - \frac{k}{2}\sum_{r\ge 4}^{k-2} E_r(q)E_{k+2-r}(q),
\end{align*}
i.e. the weight $k+2$ modular form $D_{k}E_k(q)$ is explicitly known
in this case.  
\end{remark}
\medskip

We now consider a non-linear differential equation satisfied by the Weierstrass function analogous to \eqref{g0breakformula} for rational functions.
\begin{thm}
\label{g1breakformula}
\begin{align}
\Big(\qd - P_1(x,q)\der{x} - P_1(y,q)\der{y} + P_2(x,q)+P_2(y,q)\Big)
P_2(x-y,q)
= P_2(x,q)P_2(y,q).\notag \\
\label{g1_mident1}
\end{align}
\end{thm}
\begin{prf}
Apply $\der{z}^2$ to \eqref{genus1ellthmeqn} to find
$\qd P_2(z,q) = 3P_4(z,q) -P_2^2 (z,q) - 2P_1(z,q)P_3(z,q)$.
We then find that the LHS of \eqref{g1_mident1}
is symmetric and elliptic in $x,y$ with poles only at $x=0$ and $y=0$. 
The result follows on comparing the leading terms of the
Laurant expansion in the neighbourhood of the poles on both sides of \eqref{g1_mident1}. \qed
\end{prf}
We also note the identity 
\begin{equation}
\qd E_2(q)   + 2E_2(q)P_2(x,q) + P_4(x,q)= P_2^2(x,q),
\label{p4form}
\end{equation}
obtained from  \eqref{g1_mident1} in the $x\rightarrow y$ limit.
\medskip

\subsection{Genus one $n$-point correlation functions}
\label{sec:g1_sqbracket}
In order to define and compute $n$-point correlation functions at genus one, 
Zhu \cite{Z} introduced the so-called square bracket formalism. 
This consists of a second VOA $(V,Y[\ ,\ ],\mathbf{1},\tilde{\omega})$ isomorphic to $(V,Y(\ ,\ ),\mathbf{1},\omega )$ with vertex operators
\begin{displaymath}
Y[v,z] = \sum_{n\in\mathbb{Z}}v[n]z^{-n-1} = Y(q_{z}^{L_{0}}v,q_{z}-1).
\end{displaymath}
The new conformal vector is $\wtil=\omega-\frac{C}{24}\vv$ with vertex operator
\begin{align*}
Y[\tilde{\omega},z] &= \sum_{n\in\mathbb{Z}}L[n]z^{-n-2}.
\end{align*}

The genus one partition (or $0$-point) function for $V$ is defined by 
\begin{equation}
Z_{V}^{(1)}(q)=\Tr_{V}\left(q^{L_{0}-C/24}\right)=\sum_{n\ge 0}\dim V_n q^{n-C/24}.
\label{Z1_1}
\end{equation}
For $v_{1},\ldots,v_{n}\in V$ we
define the genus one $n$-point correlation function for
parameters $z_{1},\ldots,z_{n}$ by
\begin{equation}
Z_V^{(1)}\left(v_{1},z_1;\ldots;v_{n},z_n;q\right)=
\Tr_{V}\left(Y(q_{z_{1}}^{L_{0}}v_{1},q_{z_{1}})\ldots 
Y(q_{z_{n}}^{L_{0}}v_{n},q_{z_{n}})q^{L_{0}-C/24}\right).
\label{Z1_1_npt}
\end{equation}
Zhu describes a recursion formula for expanding an $n$-point correlation function as a linear combinations of $n-1$ point correlation functions with coefficients given by elliptic functions \cite{Z}. 
Thus genus one $n$-point correlation functions are elliptic functions just as genus zero $n$-point correlation functions are rational functions.
\medskip
   
We next define the order $m$ genus one Virasoro generating function   by
\begin{align}
\notag G_m^{(1)}(z_1,\ldots,z_m;q)=Z_V^{(1)}\left(\wtil,z_1;\ldots;\wtil,z_m;q\right).
\end{align}
Analogously to Proposition~\ref{prop_Ggen0}, we then find, using Lemma~5 of \cite{MT1}, that
\begin{proposition}
\label{prop_Ggen1}
$G_m^{(1)}(z_1,\ldots,z_m;q)$ is a symmetric elliptic function in $z_1,\ldots,z_m$. 
The set of  such functions generates all genus one $n$-point correlation functions for Virasoro vacuum descendants.
\end{proposition}
Applying the Zhu recursion formula (Proposition~4.3.2, \cite{Z}) we obtain the following genus one analogue of 
Proposition~\ref{zhuredform}:
\begin{proposition}
\label{g1_zhured}
The genus one Virasoro generating function obeys the recursion relation 
\begin{align}
G_n^{(1)}(z_1,z_2,\ldots,z_n;q) &= \qd G_{n-1}^{(1)}(z_2,\ldots,z_n;q)\notag\\
&\indent+ \sum_{k=2}^nP_1(z_{1k},q)\der{z_k}G_{n-1}^{(1)}(z_2,\ldots,z_n;q)\notag\\
&\indent+ \sum_{k=2}^n2P_2(z_{1k},q)G_{n-1}^{(1)}(z_2,\ldots,z_n;q)\label{g1_zhuredthmeqn}\notag\\
&\indent+ \sum_{k=2}^n\frac{C}{2}P_4(z_{1k},q)G_{n-2}^{(1)}(z_2,\ldots,\widehat{z_k},\ldots,z_n;q),
\end{align}
where $\hat{z_k}$ denotes that $z_k$  is omitted. \qed
\end{proposition}

The RHS of \eqref{g1_zhuredthmeqn} is not manifestly symmetric in $z_1,\ldots,z_n$ as required by Proposition~\ref{prop_Ggen1}. 
We show below in Theorem~\ref{g1_maintheorem} how to express $G_n^{(1)}$ in a symmetric fashion as a sum of weights of appropriate graphs. 
In order to achieve this it is useful to define
\begin{eqnarray*}
\zh(q) &=& \eta^C(q)Z_V^{(1)}(q),\\
\Gamma_n(z_1,\ldots,z_n;q)&=&
\eta^C(q)G_n^{(1)}(z_1,z_2,\ldots,z_n;q),
\end{eqnarray*}
for Dedekind $\eta(q)$ of \eqref{intro_dedekindef}. 
Thus $\Gamma_0(q)=\zh(q)$.
Applying \eqref{deleta} we find
\begin{align*}
\Gamma_1(z_1;q) &= \left(\qd + \frac{C}{2}E_2(q)\right)\zh(q).
\end{align*}
For $n=2$ we apply \eqref{p4form} to obtain
\begin{equation*}
\Gamma_2(z_1,z_2;q)=
\left(\qd^2 + C E_2(q)\qd + \left(\frac{C}{2}E_2(q)\right)^2 + 2P_2(z_{12},q)\qd + \frac{C}{2}P_2^2(z_{12},q)\right)\zh(q).
\end{equation*}
For $n=3$ we apply \eqref{g1_mident1} to eventually obtain
\begin{align}
\Gamma_3(z_1,z_2,z_3;q)&=\Bigg(
\qd^3 + 3\frac{C}{2}E_2(q)\qd^2 + 3\left(\frac{C}{2}E_2(q)\right)^2\qd + \left(\frac{C}{2}E_2(q)\right)^3
\notag\\
&\qquad
+ 2\left(P_2(z_{12},q)+ P_2(z_{13},q) + P_2(z_{23},q)\right)\qd^2
\notag\\
&\qquad
+CE_2(q)\left(P_2(z_{12},q) + P_2(z_{13},q) + P_2(z_{23},q)\right)\qd
\notag\\
&\qquad 
+ \frac{C}{2}\left(P_2^2(z_{12},q)+ P_2^2(z_{13},q)+P_2^2(z_{23},q)\right)\qd
\notag\\
&\qquad 
+\frac{C^2}{4}E_2(q)\left(P_2^2(z_{12},q) + P_2^2(z_{13},q)+ P_2^2(z_{23},q)\right)
\notag\\
&\qquad 
+ 2\left(P_2(z_{12},q)P_2(z_{13},q) + P_2(z_{12},q)P_2(z_{23},q) + P_2(z_{13},q)P_2(z_{23},q)\right)\qd
\notag\\
&\qquad
+CP_2(z_{12},q)P_2(z_{13},q)P_2(z_{23},q)\Bigg)\zh(q).
\label{Psi3}
\end{align}
This expression depends on $C$, $E_2(q)$, $P_2(z_{ij},q)$ and $\zh(q)$ and is manifestly symmetric in $z_1,z_2,z_3$ as expected from Proposition~\ref{prop_Ggen1}. 

\subsection{Genus one Virasoro graphs}
\label{ssect:genus1}
Similarly to Section~\ref{ssect:genus0} we now develop a graphical/combinatorial approach for computing the 
order $n$ genus one Virasoro generating function. 
We define an \textit{Order $n$ Genus One Virasoro Graph} to be directed a graph $\mathcal{G}$ with $n$ vertices labelled by $z_1,\ldots,z_n$. 
Each $z_{i}$--vertex has degree $\deg(z_i)=0,1$ or $2$.
The degree $1$ vertices can have either unit indegree or outdegree 
whereas the
degree $2$ vertices have both unit indegree and outdegree. 
Thus, in this case, the connected subgraphs of $\mathcal{G}$ consist of $r$-cycles for $r\ge 1$ vertices and necklaces with degree 1 end--vertices and possible intermediate degree 2 vertices. 
We regard a single degree 0 vertex as being a degenerate necklace.

We define a weight $W(\mathcal{G})$ on $\mathcal{G}$ as follows. For each directed edge $\mc{E}_{ij}$ define an edge weight
\begin{align}
W(\mc{E}_{ij}) =W_{ij}= \left\{
  \begin{array}{ll}
    E_2(q) & i=j,\\
    P_2(z_{ij},q) & i\neq j.
  \end{array}\right.
  \label{Wij}
\end{align}
Let $K$ be the number of cycle graphs in $\mathcal{G}$ and let $M$ be the number of necklaces and define
\begin{equation}
W(\mathcal{G}) = \left(\frac{C}{2}\right)^K \prod_{\mathcal{E}}W(\mathcal{E}).\, \qd^{M}\zh(q),\label{g1_weightofgraph}
\end{equation}
where the product ranges over all the edges of \mc{G}.
Thus the weight depends on $C$, $E_2(q)$, $P_2(z_{ij},q)$ and $\zh(q)$. 

In the Appendix we describe a generating function for counting the number
$p_{n,K,M}$ of inequivalent order $n$ Virasoro graphs containing $K$ cycles and $M$ necklaces determined by the polynomial  
\begin{equation}
p_{n}(\alpha,\beta)=\sum_{K\ge 0,M\ge 0}p_{n,K,M}\alpha^M \beta^K=(-1)^n n! \sum_{i=0}^{n} \frac{(-\alpha)^i}{i!}\binom{-\beta-i}{n-i},
\label{pnKM}
\end{equation}
for necklace and cycle counting parameters $\alpha$ and $\beta$ respectively.
Thus for $n=1$ we find $p_1(\alpha,\beta) = \alpha+\beta$ corresponding to the  two inequivalent graphs
\begin{equation*}
\mathcal{G}_1=\xy (0,0)*{\cir<3pt>{}}="a"*+!R{z_1};\endxy\, ,
\qquad\mathcal{G}_2 = \xy (0,0)*[o]=<0.5pt>+{\cir<3pt>{}}="b"*+!L{z_1};
\ar@(ul,dl) "b";"b";
\endxy,
\end{equation*}
i.e. a degenerate necklace and a 1-cycle.
Clearly the sum of their weights is
\begin{equation*} 
W(\mathcal{G}_1) + W(\mathcal{G}_2) = 
\qd\zh(q) + \frac{C}{2}E_2(q)\zh(q) 
= \Gamma_1(z_1;q).
\end{equation*}
For $n=2$ we have $p_2(\alpha,\beta) = \alpha^2+2\alpha\beta+\beta^2+2\alpha+\beta$ for the 7 inequivalent graphs:
\begin{eqnarray*}
&\xy
(0,0)*[o]=<0.4pt>+{\cir<3pt>{}}="a"*+!R{z_1\,}; (10,0)*[o]=<0.4pt>+{\cir<3pt>{}}="b"*+!L{\,z_2};\endxy,
\qquad
\xy
(0,0)*[o]=<0.4pt>+{\cir<3pt>{}}="a"*+!R{z_1\,}; (10,0)*[o]=<0.4pt>+{\cir<3pt>{}}="b"*+!L{\,z_2}; \ar@(ru,lu) "a";"a";\endxy,\qquad
\xy
(0,0)*[o]=<0.4pt>+{\cir<3pt>{}}="a"*+!R{z_1\,}; (10,0)*[o]=<0.4pt>+{\cir<3pt>{}}="b"*+!L{\,z_2}; \ar@(ru,lu) "b";"b";\endxy,\qquad
\xy
(0,0)*[o]=<0.4pt>+{\cir<3pt>{}}="a"*+!R{z_1\,}; (10,0)*[o]=<0.4pt>+{\cir<3pt>{}}="b"*+!L{\,z_2}; \ar@(ru,lu) "a";"a";\ar@(ru,lu) "b";"b";\endxy,&\\
&\xy
(0,0)*[o]=<0.4pt>+{\cir<3pt>{}}="a"*+!R{z_1\,}; (10,0)*[o]=<0.4pt>+{\cir<3pt>{}}="b"*+!L{\,z_2}; \ar "a";"b";\endxy,\qquad
\xy
(0,0)*[o]=<0.4pt>+{\cir<3pt>{}}="a"*+!R{z_1\,}; (10,0)*[o]=<0.4pt>+{\cir<3pt>{}}="b"*+!L{\,z_2}; \ar "b";"a";\endxy,\qquad
\xy
(0,0)*[o]=<0.4pt>+{\cir<3pt>{}}="a"*+!R{z_1\,}; (10,0)*[o]=<0.4pt>+{\cir<3pt>{}}="b"*+!L{\,z_2}; \ar@/^/ "b";"a";\ar@/^/ "a";"b";\endxy,&
\end{eqnarray*}
whose weights  sum to
\begin{align*}
\sum_{\mathcal{G}}W(\mathcal{G})&=\qd^2\zh(q) + \frac{C}{2}E_2(q)\qd\zh(q) + \frac{C}{2}E_2(q)\qd\zh(q) + \frac{C^2}{4}E_2(q)^2\zh(q)\\
&\indent+ P_2(z_{12},q)\qd\zh(q) + P_2(z_{21},q)\qd\zh(q) + \frac{C}{2}P_2^2(z_{12},q)\zh(q)\\&=\Gamma_2(z_1,z_2;q).
\end{align*}
For $n=3$ one finds 
\begin{equation*}
p_3(\alpha,\beta) = \alpha^3+3\alpha^2\beta+6\alpha^2+3\alpha\beta^2+9\alpha\beta+6\alpha+\beta^3+3\beta^2+2\beta,
\label{eq:}
\end{equation*} 
describing 34 inequivalent genus one Virasoro graphs whose weights sum to $\Gamma_3(z_1,z_2,z_3;q)$.
These examples illustrate the general result:

\begin{thm}\label{g1_maintheorem}
$\Gamma_{n}(z_1,\ldots,z_n;q)$ is the sum of weights of all inequivalent order $n$ genus one Virasoro graphs.
\end{thm}
\begin{prf}[Proof]
We prove the result by induction in $n$. We have already shown the result holds for $n=1$ and $n=2$. For $n\ge 3$ we reinterpret the Zhu recursion formula~\eqref{g1_zhuredthmeqn} in terms of a construction of order $n$ genus one Virasoro graphs in terms of order $n-1$ and $n-2$  graphs.

Let $\mathcal{G}^{n}$ denote an order $n$ Virasoro graph.
Based on the nature of the $z_{1}$--vertex, this must belong to one of the following five types:
\begin{enumerate}[Type I.]
\item $\deg(z_{1})=0$:  
\[
\xy(0,0)*[o]=<0.4pt>+{\cir<3pt>{}}="a"*+!R{z_1\,};\endxy
\ \cdots
\]
\item $\deg(z_{1})=1$:
\[
\xy(0,0)*[o]=<0.4pt>+{\cir<3pt>{}}="a"*+!R{z_1\,};(15,0)*[o]=<0.4pt>+{\cir<3pt>{}}="c"*+!L{\,z_a};\ar "c";"a";\endxy 
\cdots\
\mbox{  or }
\xy(0,0)*[o]=<0.4pt>+{\cir<3pt>{}}="a"*+!R{z_1\,};(15,0)*[o]=<0.4pt>+{\cir<3pt>{}}="c"*+!L{\,z_a};\ar "a";"c";\endxy
\cdots
\]
\item $\deg(z_{1})=2$ where the $z_{1}$--vertex forms a 1-cycle:
\[
\xy(0,0)*[o]=<0.4pt>+{\cir<3pt>{}}="a"*+!R{z_1\,};\ar@(dr,ur) "a";"a";\endxy
\cdots
\]
\item $\deg(z_{1})=2$ where the $z_{1}$--vertex is an element of a  $2$-cycle: 
\[
\xy
(0,0)*[o]=<0.4pt>+{\cir<3pt>{}}="a"*+!R{z_1}; (10,0)*[o]=<0.4pt>+{\cir<3pt>{}}="b"*+!L{z_a}; \ar@/^/ "b";"a";\ar@/^/ "a";"b";\endxy
\cdots
\]
\item $\deg(z_{1})=2$ where the $z_{1}$--vertex is an element of a  necklace or an $r$-cycle for $r\ge 3$: 
\[
\xy
(0,0)*[o]=<0.4pt>+{\cir<3pt>{}}="a"*+!D{z_1};
(-15,0)*[o]=<0.4pt>+{\cir<3pt>{}}="b"*+!R{z_a\,};
(15,0)*[o]=<0.4pt>+{\cir<3pt>{}}="c"*+!L{\,z_b};
(-25,0)*[o]=<6pt>+{\cdots}="d";
(25,0)*[o]=<6pt>+{\cdots}="e";
\ar "a";"b"; \ar "a";"c"
\endxy
\]
\end{enumerate}
Applying the Zhu reduction formula~\eqref{g1_zhuredthmeqn} we have
\begin{align}
\Gamma_{n}(z_1,\ldots,z_n) &= \eta^{C}\qd\left(\eta^{-C}\Gamma_{n-1}(z_2,\ldots,z_n)\right)\label{Zhu1}\\
&\indent+\sum_{k=2}^n P_1(z_{1k})\der{z_{k}}\Gamma_{n-1}(z_2,\ldots,z_n)\label{Zhu2}\\
&\indent+\sum_{k=2}^n 2P_2(z_{1k}) \Gamma_{n-1}(z_2,\ldots,z_n)\label{Zhu3}\\
&\indent+\sum_{k=2}^n \frac{C}{2}P_4(z_{1k}) \Gamma_{n-2}(z_2,\ldots,\widehat{z_k},\ldots,z_n)\label{Zhu4}.
\end{align}
Here, and below, we suppress the $q$-dependence as a convenient abbreviation. 

We now show how the different parts of Zhu's reduction formula relate to the Virasoro graph weights by using induction to express  both $\Gamma\ub{n-1}{1}$ and $\Gamma\ub{n-2}{1}$ as the sum of weights of all order $n-1$ and $n-2$ Virasoro graphs, respectively. 
Let $\mathcal{G}^{n-1}$ denote an order $n-1$ Virasoro graphs labelled by $z_{2},\ldots, z_{n}$.
Then by induction we have
\begin{equation*}
\Gamma_{n-1} = \Gamma_{n-1}(z_{2},\ldots, z_{n})
=\sum_{\mc{G}^{n-1}}W(\mc{G}^{n-1}),
\end{equation*}
for
\begin{equation*}
W(\mathcal{G}^{n-1}) = (\frac{C}{2})^K\prod_{\mathcal{E}}W(\mathcal{E})\qd^{M}\zh,
\end{equation*}
where $\mathcal{G}^{n-1}$ consists of $M$ necklaces and $K$ cycles with edges $\{\mathcal{E}\}$. 

Consider the first contribution \eqref{Zhu1} arising from Zhu reduction. Recalling \eqref{deleta} then $\qd(\eta^{-C}\zh) = \eta^{-C}(\frac{C}{2}E_2 + \qd)\zh$ so that
\begin{align}
\eta^{C}\qd\left(\eta^{-C}W\left(\mathcal{G}^{n-1}\right)\right) = 
&\left(\frac{C}{2}\right)^K\left(\prod_{\mathcal{E}}W\left(\mathcal{E}\right)\right)\qd^{M+1}\zh\label{g1_pfmain5}\\
&+E_2\left(\frac{C}{2}\right)^{K+1}\left(\prod_{\mathcal{E}}W\left(\mathcal{E}\right)\right)\qd^{M}\zh\label{g1_pfmain6}\\
&+\left(\frac{C}{2}\right)^K\left(\qd\prod_{\mathcal{E}}W\left(\mathcal{E}\right)\right)\qd^{M}\zh.
\label{g1_pfmain7}
\end{align}
The first term \eqref{g1_pfmain5} is the weight of the order $n$ Virasoro graph $\mc{G}_{\mathrm{I}}^{n}$ of Type~I
consisting of a $z_{1}$-labelled degenerate necklace together with $\mathcal{G}^{n-1}$.
Clearly, all inequivalent Type~I graphs can arise in this way.
The second term \eqref{g1_pfmain6} is the weight of the order $n$ Virasoro graph $\mc{G}_{\mathrm{III}}^{n}$ of Type~III
consisting of a $z_{1}$-labelled 1-cycle together with $\mathcal{G}^{n-1}$.
Likewise, all inequivalent Type~III graphs arise in this way.
Thus summing \eqref{g1_pfmain5} and \eqref{g1_pfmain6} over all $\mathcal{G}^{n-1}$ we find 
\begin{equation}
\sum_{\mc{G}^{n-1}} \left(\frac{C}{2}\right)^K
\prod_{\mathcal{E}}W\left(\mathcal{E}\right)
\left(
\qd^{M+1}
+\frac{C}{2}E_2\qd^{M}
\right)
\zh
= 
\sum_{\mc{G}_{\mathrm{I}}^{n}} W(\mc{G}_{\mathrm{I}}^{n}) + \sum_{\mc{G}_{\mathrm{III}}^{n}} W(\mc{G}_{\mathrm{III}}^{n}),
\label{sumI_III}
\end{equation}
summing over all inequivalent order $n$ graphs of Type~I and III.
The remaining term \eqref{g1_pfmain7} contributes
\begin{eqnarray}
\left(\frac{C}{2}\right)^K\left(\qd\prod_{\mathcal{E}}W\left(\mathcal{E}\right)\right)\qd^{M}\zh
&=&
W(\mathcal{G}^{n-1})
\sum_{\mathcal{E}}\qd \log W(\mathcal{E}).
\label{qdelW}
\end{eqnarray}
In a similar fashion, \eqref{Zhu2}  of the Zhu reduction formula gives a $\mathcal{G}^{n-1}$ contribution of 
\begin{equation}
W(\mathcal{G}^{n-1})\sum_{\mathcal{E}}
\sum_{k=2}^n 
P_1(z_{1k})\der{z_{k}}\log W(\mathcal{E}).
\label{zdelW}
\end{equation}
Likewise, \eqref{Zhu3} gives a $\mathcal{G}^{n-1}$ contribution of 
\begin{equation}
W(\mathcal{G}^{n-1})\sum_{k=2}^n 2 P_2(z_{1k})=
W(\mathcal{G}^{n-1})
\left(
\sum_{\mathcal{E}_{ab}}\big(P_2(z_{1a})+P_2(z_{1b})\big)+
\sum_{c,\,\deg(z_c)=1} P_2(z_{1c})
\right),
\label{P2W}
\end{equation}
where the first sum is taken over all $\mathcal{G}^{n-1}$ edges $\mathcal{E}_{ab}$ and
 the second sum over all $z_c$--vertices of degree 1 arising from  necklace end-vertices. 
Let 
\begin{eqnarray}
S(\mathcal{G}^{n-1})&=&W(\mathcal{G}^{n-1})
\Big( 
\sum_{\mathcal{E}_{ab}}
\left[
\big(
 \qd + \sum_{k=2}^n 
P_1(z_{1k})\der{z_{k}}
\big) \log W(\mathcal{E}_{ab})
+P_2(z_{1a})+P_2(z_{1b})
\right]
\notag\\
&&\quad 
+\sum_{c,\,\deg(z_c)=1} P_2(z_{1c})
\Big),
\label{SG}
\end{eqnarray}
denote the sum of \eqref{qdelW}--\eqref{P2W}. 
Consider an edge $\mathcal{E}_{ab}$ with $z_a,z_b$--vertices contributing to \eqref{SG} 
where either (A)~$z_a=z_b$ or (B)~$z_a \neq z_b$.

\textbf{Case (A).}
Here the $z_{a}$--vertex forms a 1-cycle with edge weight $W(\mathcal{E}_{aa})=E_2$. 
Hence $W(\mathcal{G}^{n-1})=\frac{C}{2}E_2 W(\mathcal{G}^{n-2})$ 
where $\mathcal{G}^{n-2}$ is an order $n-2$ graph labelled by $z_2,\ldots, \widehat{z_{a}},\ldots, z_{n}$ (with label $\widehat{z_{a}}$ deleted).
All distinct such order $n-2$ graphs arise in this way.  
The $\mathcal{E}_{aa}$ dependent part of \eqref{SG} is thus
\begin{equation}
\frac{C}{2}W(\mathcal{G}^{n-2})E_2 \left(\qd \log E_2+2P_2(z_{1a})\right).
\label{E2term}
\end{equation}
The last term \eqref{Zhu4} in the Zhu reduction formula  makes a
$\mathcal{G}^{n-2}$ contribution of $\frac{C}{2}W(\mathcal{G}^{n-2})P_4(z_{1a})$.
Adding this to \eqref{E2term} and using \eqref{p4form} we obtain
\begin{eqnarray*}
W(\mathcal{G}^{n-2})
\frac{C}{2}
\big(\qd E_2  +2E_2P_2(z_{1a})+P_4(z_{1a})\big)
&=&
W(\mathcal{G}^{n-2})
\frac{C}{2}P_2(z_{1a})^2\\
&=&
W(\mathcal{G}^{n}_{\mathrm{IV}}),
\end{eqnarray*}
where $\mathcal{G}^{n}_{\mathrm{IV}}$ is the order $n$ graph of Type~IV formed by the union of $\mathcal{G}^{n-2}$
with the 2-cycle with vertices labelled by $z_1$ and $z_a$.
Summing over all $z_a$ and $\mathcal{G}^{n-2}$ we  obtain 
$
\sum_{\mc{G}_{\mathrm{IV}}^{n}} W(\mc{G}_{\mathrm{IV}}^{n})$.

\textbf{Case (B).}  In this case $W(\mathcal{E}_{ab})=P_2(z_{ab})$
so that the total $\mathcal{E}_{ab}$ contribution arising from the edge sum in $S(\mathcal{G}^{n-1})$ of \eqref{SG} is
\begin{eqnarray}
&&W(\mathcal{G}^{n-1})P_{2}(z_{ab})^{-1}
\Big(
\qd 
+P_1(z_{1a})\der{z_{a}}
+P_1(z_{1b})\der{z_{b}}
+P_2(z_{1a})+P_2(z_{1b})
\Big)P_{2}(z_{ab})\notag\\ &&
=W(\mathcal{G}^{n-1})\frac{P_2(z_{1a})P_2(z_{1b})}{P_{2}(z_{ab})},
\label{P2abeqn}
\end{eqnarray}
using Theorem~\ref{g1breakformula} with $x=-z_{1a},y=-z_{1b}$. 
But the RHS of \eqref{P2abeqn} is the weight of the order $n$ 
graph $\mathcal{G}^{n}_{\mathrm{V}}$ of type~V constructed by 
inserting $z_{1}$ between the $z_{a}$--vertex and the $z_{b}$--vertex: 
\[
\cdots\xy
(0,0)*[o]=<0.4pt>+{\cir<3pt>{}}="a"*+!R{z_a\,}; 
(15,0)*[o]=<0.4pt>+{\cir<3pt>{}}="d"*+!D{z_1}; 
(30,0)*[o]=<0.4pt>+{\cir<3pt>{}}="b"*+!L{\,z_b};\ar "b";"d";
\ar "d";"a";\endxy
\cdots
\]
Summing over inequivalent $\mathcal{G}^{n-1}$ we obtain 
$
\sum_{\mc{G}_{\mathrm{V}}^{n}} W(\mc{G}_{\mathrm{V}}^{n}),
$ 
summing over all inequivalent graphs of type~V. 

Lastly, if  $\deg{z_{a}}=1$ then the $z_{a}$--vertex is a necklace
end--vertex and we obtain an additional contribution
from the last term in \eqref{SG} of
\begin{equation*}
W(\mathcal{G}^{n-1})P_2(z_{1a})=W(\mathcal{G}^{n}_{\mathrm{II}}),
\end{equation*} 
where $\mathcal{G}^{n}_{\mathrm{II}}$ is the order $n$ graph of Type~II 
found by joining the $z_{1}$--vertex to the $z_{a}$--vertex. Summing over inequivalent $\mathcal{G}^{n-1}$ and all such necklace end $z_{a}$--vertices results in 
$
\sum_{\mc{G}_{\mathrm{II}}^{n}} W(\mc{G}_{\mathrm{II}}^{n})
$. 
Thus altogether we have shown that 
\begin{equation*}
\Gamma_{n}(z_1,\ldots,z_n;q)=\sum_{\mc{G}^{n}} W(\mc{G}^{n}),
\end{equation*}
summing over all inequivalent order $n$ Virasoro graphs of type~I--V and hence the theorem holds.
\qed
\end{prf}
\begin{remark}
Theorem~\ref{g1_maintheorem} can be readily generalized for any $V$-module $M$ by replacing $\Theta_V(q)$ by 
$\Theta_M(q)=\eta(q)^{C}\Tr_M \left(q^{L(0)-C/24}\right)$ throughout. 
The remaining dependence on $C$, $E_2(q)$ and $P_2(z_{ij},q)$ remains unchanged.
\end{remark}

\medskip
\subsection{Partial permutations}
\label{sec:pperm}
Let $\Psi_{n}$ denote the set of partial permutations of the label set $\Phi=\lbrace 1,\ldots,n\rbrace$  
i.e. the injective partial mappings (including the empty map) from $\{1,\ldots ,n\}$ to itself.  
The set of order $n$ genus one Virasoro graphs is in 1-1 correspondence with $\Psi_{n}$ where the directed edges $\{\mathcal{E}_{ij}\}$ of $\mathcal{G}$ correspond to $\psi\in \Psi_n$ with $\psi(i)=j$. 
 
It is natural to introduce the \emph{Partial Permanent} of an $n\times n$ matrix $A=(a_{ij})$ indexed by $i,j\in \{1,\ldots ,n\}$ as follows 
\begin{equation}
\pperm\, A= \sum\limits_{\psi\in\Gamma }\prod\limits_{i \in\, \dom \psi} a_{i\psi(i)},
\label{pperm}
\end{equation}
where $\dom \psi$ denotes the domain of the partial map $\psi$ and
with unit contribution for the empty map. 
We also define the \emph{$\alpha,\beta$-Extended Partial Permanent} by 
\begin{equation}
\pperm_{\alpha\beta} A= \sum\limits_{\psi\in\Psi_{n}}
\alpha^{N(\psi)}\beta^{K(\psi)}
\prod\limits_{i \in \dom \psi} a_{i\psi(i)}
,
\label{ppermBTP}
\end{equation}
where $N(\psi)$ is the number of necklaces and $K(\psi)$ is the number of cycles determined by $\psi$. 

We now utilize partial permutations to re-express the order $n$ genus one Virasoro generating function. 
Towards this end,  we define a linear map $\mathcal{L}_{\alpha}$ from $\mathbb{C}[\alpha]$, the vector space of complex coefficient polynomials in $\alpha$,  to the complex vector space spanned by $\qd^k\zh(q)$ as follows
\begin{align}
   \mathcal{L}_{\alpha}(\sum_{k=0}^n c_k \alpha^k) &= \sum_{k=0}^n
   c_k \qd^k\zh(q), \quad c_k\in \mathbb{C}.
   \label{Lmap}
\end{align}
Then Theorem~\ref{g1_maintheorem} can be restated as
\begin{thm}
\label{g1mainthmperm}
The order $n$ genus one Virasoro generating function is given by
\begin{displaymath}
\Gamma_n(z_1,\ldots,z_n;q) = \mathcal{L}_{\alpha}\left(\pperm_{\alpha \frac{C}{2}}W\right),\qed
\end{displaymath}
with $W=(W_{ij})$ for edge weight $W_{ij}$ of \eqref{Wij}. 
\end{thm}

\medskip
\subsection{Further combinatorics}
\label{sec:g1_combin} 
Lastly, we briefly consider an alternative way of 
computing $\Gamma_{n}(z_1,\ldots,z_n;q)$ in terms of a sum over the  permutations of the label set $\Phi = \lbrace 1,\ldots,n\rbrace$
by means of an alternative graphical approach. 
Define a truncation operation $\mathcal{T}_{\rho}$ on the vector space of polynomials $\mathbb{C}[\rho]$ by  
\begin{equation*}
\mathcal{T}_{\rho}(\sum_{k=0}^n c_k \rho^k) = c_0 + c_1 \rho,\quad c_k\in \mathbb{C}.
\end{equation*}
We also define a linear map $\mathcal{L}_{\rho}$ as in \eqref{Lmap}.

We next define an \textit{Order $n$ Genus One Permutation Virasoro Graph}
to be a directed graph $\mathcal{G}$ with $n$ vertices labelled by $z_1,\ldots,z_n$ where each vertex has degree 2 with unit indegree and outdegree. 
The set of such graphs is in one to one correspondence with the permutation group $\Sigma_n$ on $n$ symbols.  
We define a weight $\chi(\mathcal{G})$ on $\mathcal{G}$ as follows. For each directed edge $\mc{E}_{ij}$ define a weight
\begin{align*}
\chi_{ij} = \left\{
  \begin{array}{ll}
    \frac{2}{C}\rho  +E_2(q) & i=j,\\
    \frac{2}{C}\rho + P_2(z_{ij},q) & i\neq j.
  \end{array}\right.
\end{align*}
Then for a cycle subgraph $\mathcal{G}_{\sigma}$ with edges $\mathcal{E}^\sigma$ we define:
\begin{displaymath}
\chi(\mathcal{G}_{\sigma}) = \frac{C}{2} \mathcal{T}_{\rho}(\prod_{\mathcal{E}^\sigma_{ij}}\chi_{ij}).
\end{displaymath} 
Each $\sigma$ corresponds to permutation group cycle.
Finally we define the weight of $\mc{G}$ by
\begin{align*}
\chi(\mathcal{G}) =
\mathcal{L}_{\rho}\left(\prod_{\sigma}\chi(\mathcal{G}_{\sigma})\right).
\end{align*}
For example, there are $3!$ order $3$ permutation graphs
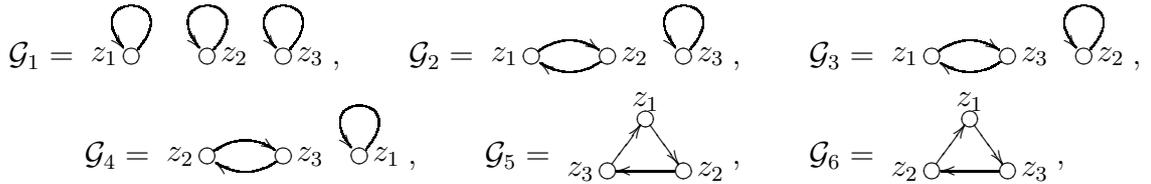
\begin{eqnarray*}
&
\mathcal{G}_1=
\xy
(0,0)*[o]=<0.4pt>+{\cir<3pt>{}}="a"*+!R{z_1\,}; (10,0)*[o]=<0.4pt>+{\cir<3pt>{}}="b"*+!L{\,z_2}; 
(20,0)*[o]=<0.4pt>+{\cir<3pt>{}}="c"*+!L{\,z_3}; 
\ar@(ru,lu) "a";"a";
\ar@(ru,lu) "b";"b";
\ar@(ru,lu) "c";"c";
\endxy,\qquad
\mathcal{G}_2= 
\xy
(0,0)*[o]=<0.4pt>+{\cir<3pt>{}}="a"*+!R{z_1\,}; (10,0)*[o]=<0.4pt>+{\cir<3pt>{}}="b"*+!L{\,z_2}; 
(20,0)*[o]=<0.4pt>+{\cir<3pt>{}}="c"*+!L{\,z_3};
\ar@/^/ "a";"b";
\ar@/^/ "b";"a";
\ar@(ru,lu) "c";"c";
\endxy,\qquad 
\mathcal{G}_3=
\xy
(0,0)*[o]=<0.4pt>+{\cir<3pt>{}}="a"*+!R{z_1\,}; (10,0)*[o]=<0.4pt>+{\cir<3pt>{}}="b"*+!L{\,z_3}; 
(20,0)*[o]=<0.4pt>+{\cir<3pt>{}}="c"*+!L{\,z_2};
\ar@/^/ "b";"a";
\ar@/^/ "a";"b";
\ar@(ru,lu) "c";"c";
\endxy,
&
\\ 
&
\mathcal{G}_4=
\xy
(0,0)*[o]=<0.4pt>+{\cir<3pt>{}}="a"*+!R{z_2\,}; (10,0)*[o]=<0.4pt>+{\cir<3pt>{}}="b"*+!L{\,z_3}; 
(20,0)*[o]=<0.4pt>+{\cir<3pt>{}}="c"*+!L{\,z_1};
\ar@/^/ "b";"a";
\ar@/^/ "a";"b";
\ar@(ru,lu) "c";"c";
\endxy,\qquad 
\mathcal{G}_5=
\xy (0,-2)*[o]=<0.5pt>+{\cir<3pt>{}}="a"*+!R{z_3\,}; (10,-2)*[o]=<0.5pt>+{\cir<3pt>{}}="b"*+!L{\,z_2}; (5,5)*[o]=<0.5pt>+{\cir<3pt>{}}="c"*+!D{z_1};
\ar "c";"a";\ar "b";"c";\ar "a";"b";\endxy,\qquad
\mathcal{G}_6=
\xy (0,-2)*[o]=<0.5pt>+{\cir<3pt>{}}="a"*+!R{z_2\,}; (10,-2)*[o]=<0.5pt>+{\cir<3pt>{}}="b"*+!L{\,z_3}; (5,5)*[o]=<0.5pt>+{\cir<3pt>{}}="c"*+!D{z_1};
\ar "c";"a";\ar "b";"c";\ar "a";"b";
\endxy,
&
\end{eqnarray*} 
in contrast to the 34 order 3 graphs arising in the previous section. 
The sum of the weights for the displayed 6 graphs results in $\Gamma_3(z_1,z_2,z_3;q)$ of \eqref{Psi3}. 
In particular,
$\chi(\mathcal{G}_1)$ gives the first line, 
$\chi(\mathcal{G}_2)+\chi(\mathcal{G}_3)+\chi(\mathcal{G}_4)$ gives lines 2--5
and $\chi(\mathcal{G}_5)+\chi(\mathcal{G}_6)$ gives lines 6--7 of \eqref{Psi3}.
This illustrates the general result:
\begin{thm}\label{genus1theoremvir2}
$\Gamma\ub{n}{1}(z_1,\ldots,z_n;q)$ is the sum of weights of all inequivalent order $n$ genus one permutation Virasoro graphs. \qed
\end{thm}
\medskip
The proof proceeds along similar lines to Theorem~\ref{g1_maintheorem}.
Using the correspondence with $\Sigma_n$ we can also restate this result as a sum over permutations:
\begin{align}
\Gamma_n(z_1,\ldots,z_n;q) =
\mathcal{L}_{\rho}\left(
\sum_{\pi\in \Sigma_n}
\left(\frac{C}{2}\right)^K
\prod_{\sigma_{k}}\mathcal{T}_{\rho}\left(\prod_{i\in\sigma_{k}}\chi_{i\sigma(i)}\right)
\right)
\label{g1_combformulanpt},
\end{align}
where $\pi$ has cycle decomposition  $\pi=\sigma_{1}\ldots \sigma_{K}$ with  cycles  $\sigma_{k}=(\ldots i\ldots)$.

\appendix

\section{Counting Graphs}
\subsection{Counting genus zero Virasoro graphs}
Let $d_{n,K}$ denote the number of order $n$ genus zero Virasoro graphs with $K$ cycles as introduced in Section~\ref{ssect:genus0}.  
Define a polynomial $d_n(\beta)=\sum_{K\ge 1}d_{n,K} \beta^K$ with $d_0(\beta)=1$. 
Recall that each order $n$ genus zero graph can be expressed as a derangement on $n$ symbols. 
Hence $d_n(1)$ is the number of derangements on $n$ symbols. 
Define an exponential generating function for genus zero Virasoro graphs
\begin{equation*}
D(\beta,z) =\sum_{n\ge 0}\frac{z^n}{n!} d_n(\beta).
\end{equation*}
We then find
\begin{proposition}
\begin{equation}
D(\beta,z) = \left(\frac{e^{-z}}{1-z}\right)^\beta.
\label{dbetaz}
\end{equation}
\end{proposition}
\begin{prf}
Define a weight function $w_{\beta,z}$ on each order $n$ graph 
with $K$ cycles by
\begin{equation*}
w_{\beta,z}(\mathcal{G})=\beta^K z^n.
\end{equation*}
The exponential generating function is thus
\begin{equation*}
D(\beta,z) =\sum_{\mathcal{G}}\frac{1}{n!}w_{\beta,z}(\mathcal{G}),
\end{equation*}
where $n$ denotes the order of $\mathcal{G}$. 
Each genus zero Virasoro graph $\mathcal{G}$ is the 
labelled product of $r$-cycles for $r\ge 2$.
But $w_{\beta,z}$ is multiplicative with respect to this cycle decomposition with $w_{\beta,z}(\mathcal{G}_{\sigma})=\beta z^r$ for any $r$-cycle. 
Thus we find e.g. \cite{FS} 
\begin{equation*}
D(\beta,z)= \exp(C(\beta,z)),
\end{equation*}
where $C(\beta,z)=\sum_{\mathcal{G}_{\sigma}}\frac{1}{r!}w_{\beta,z}(\mathcal{G}_{\sigma})$ is the exponential generating function for all $r$-cycles for $r\ge 2$. Since there are $(r-1)!$ inequivalent $r$-cycles this is given by
\begin{align*}
C(\beta,z) &= \sum_{r\ge 2}\frac{1}{r!}\beta z^r (r-1)!
=\beta(-z-\log(1-z)).
\end{align*}
Therefore the result follows. \qed
\end{prf}

It is straightforward to expand \eqref{dbetaz} to obtain 
\eqref{dnbeta}
and the following recursion formula
\begin{align*}
d_n(\beta) = (n-1)&(d_{n-1}(\beta)+\beta d_{n-2}(\beta)),\quad d_{-1}(\beta) = 0,\, d_0(\beta)=1.
\end{align*}
These agree with standard results for derangements in the case $\beta=1$ e.g. \cite{FS}. 
\medskip

\subsection{Counting genus one Virasoro graphs}
Let $p_{n,K,M}$ denote the number of order $n$ genus one Virasoro graphs with $K$ cycles and $M$ necklaces as introduced in Section~\ref{ssect:genus1}.  
Define a polynomial $p_n(\alpha,\beta)
=\sum_{K\ge 0,M\ge 0}p_{n,K,M}\alpha^M \beta^K$ with $p_0(\alpha,\beta)=1$.  
Define the exponential generating function for genus one Virasoro graphs
\begin{equation*}
P(\alpha,\beta,z) =\sum_{n\ge 0}\frac{z^n}{n!} p_n(\alpha,\beta).
\end{equation*}
We then find
\begin{proposition}
\begin{equation}
P(\alpha,\beta,z) = \frac{\exp( \frac{\alpha z}{1-z})}{(1-z)^\beta}.
\label{pbetaz}
\end{equation}
\end{proposition}
\begin{prf}
Here we define a weight function $w_{\alpha,\beta,z}$ on each order $n$ graph 
with $K$ cycles and $M$ necklaces by
\begin{equation*}
w_{\alpha,\beta,z}(\mathcal{G})=\alpha^M \beta^K z^n,
\end{equation*}
so that
\begin{equation*}
P(\alpha,\beta,z) =\sum_{\mathcal{G}}\frac{1}{n!}w_{\alpha,\beta,z}(\mathcal{G}),
\end{equation*}
summing over all genus one Virasoro graphs $\mathcal{G}$ for all order $n$. 
But the weight function is multiplicative with respect to the necklace $\mathcal{G}_n$
and cycle $\mathcal{G}_{\sigma}$ decomposition so that here we find
e.g. \cite{FS}
\begin{eqnarray*}
&P(\alpha,\beta,z)=  \exp(N(\alpha,z)+C(\beta,z)),
\end{eqnarray*}
where $N(\alpha,z)=\sum_{\mathcal{G}_N}\frac{1}{r!}w_{\alpha,\beta,z}(\mathcal{G}_N)$ 
is the exponential generating function for connected necklaces
and  $C(\beta,z)=\sum_{\mathcal{G}_{\sigma}}\frac{1}{r!}w_{\alpha,\beta,z}(\mathcal{G}_{\sigma})$ is the exponential generating function for connected $r$-cycles for $r\ge 1$. Thus, much as before 
\begin{align*}
C(\beta,z) &= \beta\sum_{r\ge 1}\frac{1}{r!} z^r(r-1)!
=-\beta\log(1-z).
\end{align*}
Similarly, $w_{\alpha,\beta,z}(\mathcal{G}_N)=\alpha z^r$ for a necklace with $r\ge 1$ vertices and since there are $r!$ inequivalent such necklaces we find
\begin{align*}
N(\alpha,z) &= \sum_{r\ge 1}\frac{1}{r!}\alpha z^r r!=\alpha\frac{z}{1-z}.
\end{align*}
Therefore the result follows. \qed 
\end{prf}
\eqref{pbetaz} implies \eqref{pnKM} and the following recursion formula
\begin{displaymath}
p_{n+1}(\alpha,\beta) = (2n+\alpha+\beta)p_n(\alpha,\beta) - n(n+\beta-1)p_{n-1}(\alpha,\beta),
\end{displaymath}
with $p_{-1}(\alpha,\beta) = 0$ and $p_0(\alpha,\beta)=1$.

\end{document}